\documentclass[11pt]{article}
\usepackage{amsmath, amsthm, amscd, amsfonts, amssymb, graphicx}
\usepackage{titles}
\usepackage{algorithm}
\usepackage{algpseudocode}
\usepackage{subfigure} 
\usepackage{caption}
\usepackage{graphics,graphicx}
\usepackage{amsmath, amssymb, amsthm}
\usepackage{mathrsfs} 
\usepackage{enumitem}

\usepackage{indentfirst}

\numberwithin{equation}{section}
\usepackage{epsfig}
\usepackage{color}
\usepackage{fancyhdr}
\usepackage{titlesec}
\usepackage{tikz}
\usepackage{afterpage}
\usetikzlibrary{arrows}
\usetikzlibrary{decorations.markings}
\tikzstyle{block}=[draw opacity=0.7,line width=1.4cm]
\tikzset{
	big black arrow/.style={
		decoration={markings,mark=at position 1 with {\arrow[scale=2.5,black]{>}}},
		postaction={decorate},
		shorten >=0.4pt},
	line/.style={draw, ->}}

\usepackage{eufrak}
\newtheorem{theorem}{Theorem}[section]

\newtheorem{definition}[theorem]{Definition}
\newtheorem{corollary}[theorem]{Corollary}

\newtheorem{remark}[theorem]{Remark}

\title{\textbf{On partitions with bounded multiplicity restrictions on even parts and with the largest part being odd}}
\date{}
\newcommand{\qPoch}[3]{(#1;#2)_{#3}}
\begin{document}

\maketitle

\vspace{-1.7cm}

\maketitle
\begin{center}
	{\bf Gaurab Bardhan$^1$ and Nipen Saikia$^{2, \ast}$}\\
	$^1$Department of Mathematics, Tyagbir Hem Baruah College,\\ Jamugurihat, Sonitpur, Assam, India.\\
	E. Mail: gaurabbardhan561@gmail.com
	\vskip2mm
	$^2$Department of Mathematics, Rajiv Gandhi
	University,\\ Rono Hills, Doimukh, Arunachal Pradesh, India.\\
	E. Mail(s): nipennak@yahoo.com\\
	$^\ast$\textit{Corresponding author}.\end{center}
\begin{abstract}
Andrews and El Bachraoui (\textit{RACSAM} \textbf{119}, \textit{Article Number 4 (2025)}) established identities relating $4$-regular partitions to partitions with distinct even parts whose largest part is odd. As an extension to their work, in this paper we define partitions of a positive integer $n$  with bounded multiplicity restriction on even parts and whose largest part is odd. We prove some identities relating our  partition functions and $2r$-regular partition function,  and show that identities due to Andrews and El Bachraoui follow as particular cases. 
\end{abstract}

\noindent{\bf Keywords and Phrases:} Integer partitions; partitions with bounded multiplicity restrictions;  generating functions;  $q$-series identities.\vskip2mm

\noindent{\bf Mathematics Subject Classifications: }11P81; 05A17; 11D09.

\section{Introduction}
Throughout this paper, let $a$ and $q$ be a complex numbers satisfying $|q|<1$. Define $q$-series \cite{1, 7}  in standard notation as 
\begin{equation}
(a;q)_0=1,\qquad
(a;q)_n=\prod_{j=0}^{n-1}(1-aq^j)\quad \text{and}\quad
(a;q)_\infty=\prod_{j=0}^{\infty}(1-aq^j).
\label{Bin1}
\end{equation}
The following elementary consequences of \eqref{Bin1} will be used in the paper without referring explicitely:
\begin{align}
(a;q)_{m+n}
=(a;q)_m(aq^m;q)_n,\quad
(a;q)_\infty
=(a;q)_n(aq^n;q)_\infty,\quad \text{and}\quad
(a;q)_\infty
=(a;q^2)_\infty(aq;q^2)_\infty,\notag
\end{align}
where $m$ and $n$ are positive integers.
A non-increasing finite sequence of positive integers $\beta_1\geq \beta_2 \geq \cdots \geq \beta_k>0$ is said to be a partition of a positive integer $n$ if $n=\sum_{i=1}^{k}\beta_i$.
The integers $\beta_i$ are called the parts of the partition.  
If $p(n)$ denotes the number of partitions of $n$, then its  generating function \cite{euler1748introductio} is given by
	$$\sum_{n=0}^{\infty}p(n)q^n=\dfrac{1}{(q;q)_\infty};\qquad p(0)=1.$$
For example, the partitions of $p(7) =15$ and the relevant partitions are given by $7, 6 + 1, 5 + 2, 5 + 1 + 1, 4 + 3, 4 + 2 + 1, 4 + 1 + 1 + 1, 3 + 3 + 1, 3 + 2 + 2, 3 + 2 + 1 + 1, 3 + 1 + 1 + 1 + 1, 2 + 2 + 2 + 1, 2 + 2 + 1 + 1 + 1, 2 + 1 + 1 + 1 + 1 + 1$, and  $1 + 1 + 1 + 1 + 1 + 1 + 1$.

For any integers $l\geq2$, an $l$-regular partition of a positive integer $n$ is a partition of $n$ wherein none of the  parts is divisible by $l$. Let $\mathcal R_l(n)$ denote the number of $l$-regular partitions of $n$, then its generating function is given by 
$$\sum_{n=0}^{\infty}\mathcal R_l(n)q^n=\dfrac{(q^l;q^l)_\infty}{(q;q)_\infty}.$$
For example, the partitions of $\mathcal R_4(7)=12$ and the relevant partitions are given by $7, 6 + 1, 5 + 2, 5 + 1 + 1, 3 + 3 + 1, 3 + 2 + 2, 3 + 2 + 1 + 1, 3 + 1 + 1 + 1 + 1, 2 + 2 + 2 + 1, 2 + 2 + 1 + 1 + 1, 2 + 1 + 1 + 1 + 1 + 1$, and  $1 + 1 + 1 + 1 + 1 + 1 + 1$.

Also, for any nonnegative integer $n$ let $ped(n)$  denote the number of partitions of $n$ with distinct even parts. Then  the  generating function of $ped(n)$ \cite{2} is given by 
\begin{equation}\label{pp}\sum_{n=0}^{\infty}ped(n)q^n
=
\dfrac{(-q^2;q^2)_\infty}{(q;q^2)_\infty}
=
\dfrac{(q^4;q^4)_\infty}{(q;q)_\infty}.\end{equation}
For example, the partitions of $ped(7) =12$ and the relevant partitions are given by $7, 6 + 1, 5 + 2, 5 + 1 + 1, 4 + 3, 4 + 2 + 1, 4 + 1 + 1 + 1, 3 + 3 + 1 , 3 + 2 + 1 + 1, 3 + 1 + 1 + 1 + 1,  2 + 1 + 1 + 1 + 1 + 1$, and  $1 + 1 + 1 + 1 + 1 + 1 + 1$.
The second equality of \eqref{pp} follows from the first equality by elementary $q$-series manipulations. Therefore, one clearly has, $\mathcal R_4(n)=ped(n).$ For the arithmetic properties and related result of $ped(n)$ see \cite{2, 4,6,11,12}  and for $\mathcal R_l(n)$ see \cite{CalkinEtAl,CuiGu,FurcyPenniston,GordonOno,HirschhornSellers}.

Recently, Andrews and El Bachraoui \cite{DE} derived some relations relating certain restricted partitions into distinct even parts with restricted $4$-regular partitions.  They defined  three partition functions,  \(\operatorname{DE1}(n)\) which counts the number of partitions of \(n\) in which even parts are distinct and the largest part is odd;  \(\operatorname{DE2}(n)\) which counts the number of partitions of \(n\) in which even parts are distinct, the largest part is odd and appears at least twice; and  \(\operatorname{DE3}(n)\) which counts the number of partitions of \(n\) in which even parts are distinct, the largest part is odd and appears exactly once. For example, \(\operatorname{DE1}(7)=7\) with the partitions given by \(7\), \(5+2\), \(5+1+1\), \(3+3+1\), \(3+2+1+1\), \(3+1+1+1+1\), \(1+1+1+1+1+1+1\); \(\operatorname{DE2}(7)=2\) with the partitions given by \(3+3+1\), \(1+1+1+1+1+1+1\); and \(\operatorname{DE3}(7)=5\) with the partitions given by \(
7\), \(5+2\), \(5+1+1\), \(3+2+1+1\), \(3+1+1+1+1.\)
The corresponding generating functions \cite{DE} are defined as follows:
$$\sum_{n=0}^{\infty}\operatorname{DE1}(n)q^{n}=\sum_{n=0}^{\infty}\dfrac{\qPoch{-q^{2}}{q^{2}}{n}q^{2n+1}}{\qPoch{q}{q^{2}}{n+1}},$$
$$\sum_{n=0}^{\infty}\operatorname{DE2}(n)q^{n}=\sum_{n=0}^{\infty}\dfrac{\qPoch{-q^{2}}{q^{2}}{n}q^{4n+2}}{\qPoch{q}{q^{2}}{n+1}},$$
 and
$$\sum_{n=0}^{\infty}\operatorname{DE3}(n)q^{n}=\sum_{n=0}^{\infty}\dfrac{\qPoch{-q^{2}}{q^{2}}{n}q^{2n+1}}{\qPoch{q}{q^{2}}{n}}.$$
They further derived the following partition relations \cite[Corollary 1, 2, 3 and 4]{DE}, respectively:\\
(i) For \(n>0\), \(\operatorname{DE1}(n)+\operatorname{DE1}(n-1)\) equals the number of the $4$-regular partitions of \(n\).\\
(ii) For \(n>0\), \(\operatorname{DE2}(n)+\operatorname{DE2}(n-3)\) equals the number of \(4\)-regular partitions of \(n\) into parts each \(>1\).\\
(iii) For \(n>1\), \(\operatorname{DE3}(n+2)+\operatorname{DE3}(n-1)\) equals the number of \(4\)-regular partitions of \(n\).\\
(iv) For \(n>1\), \(\operatorname{DE3}(n+2)+\operatorname{DE3}(n-1)=\operatorname{DE1}(n)+\operatorname{DE1}(n-1)\).

As an extension to the work done in \cite{DE}, in this paper we define partitions of a positive integer $n$  with bounded multiplicity restrictions on even parts and with the largest part being odd. We establish  their generating functions and prove some  identities relating the newly defined partition  functions and $2r$-regular partition function. Moreover, we  establish  \cite[Corollary 1, 2, 3 and 4]{DE} as particular cases of our results. 

The paper is organized as follows. In Section 2, we introduce our partition functions and  establish their generating functions.  In Section 3, we prove some identities relating  our partition fucntions  and $2r$-regular  partition function.   In Remark \ref{rem} we show that  \cite[Corollary 1, 2, 3 and 4]{DE} follow as  particular cases of our results.

\section{Partition functions  and their generating functions}
\begin{definition}
	Let $h\geq0$ and  $r\geq2$  be integers. Let 
	$\mathcal R_{2r}(n,h)$ to be the number of partitions of nonnegative integer $n$ into parts not divisible by $2r$ and having exactly $h$ odd parts, counted with multiplicity. 
\end{definition}For example,  $\mathcal R_4(5,1)=3$ with the relevant partitions  $5, 3+2, 2+2+1$, and  $\mathcal R_6(8,0)=4$ with the relevant partitions $8, 4+4, 4+2+2, 2+2+2+2$.

\begin{theorem}
	For any integers $h\geq0$, $r\geq2$, and $|xq|<1$, we have
	\begin{equation}\label{RGF}
		\sum_{n=0}^{\infty}\sum_{h=0}^{\infty}
		\mathcal R_{2r}(n,h)x^h q^n
		=
		\dfrac{(q^{2r};q^{2r})_{\infty}}
		{(xq;q^2)_{\infty}(q^2;q^2)_{\infty}}.
	\end{equation}
\end{theorem}
\begin{proof}Note that,  no part is divisible by $2r$ in $2r$-regular partitions of $n$.
	Also, since every multiple of $2r$ is even, all odd parts
	$1,3,5,\ldots$ are allowed without any restriction. If 
	$x$ counts the number of odd parts in the partitions, their contribution to $\mathcal R_{2r}(n,h)$ is given by 
	\begin{equation}\label{RO1}
		\prod_{j=0}^{\infty}\dfrac{1}{1-xq^{2j+1}}
		=
		\dfrac{1}{(xq;q^2)_{\infty}}.
	\end{equation}
	Also, since the  allowed even parts are $2j$ with $r\nmid j$,  their 
	contribution  to $\mathcal R_{2r}(n,h)$ is given by 
	\begin{equation}\label{RE1}
		\prod_{\substack{j\geq1\\ r\nmid j}}
		\dfrac{1}{1-q^{2j}}
		=
		\dfrac{\displaystyle\prod_{j=1}^{\infty}(1-q^{2rj})}
		{\displaystyle\prod_{j=1}^{\infty}(1-q^{2j})}
		=
		\dfrac{(q^{2r};q^{2r})_{\infty}}
		{(q^2;q^2)_{\infty}}.
	\end{equation}
	Multiplying \eqref{RO1} and \eqref{RE1}, we obtain
	$$\sum_{n=0}^{\infty}\sum_{h=0}^{\infty}
	\mathcal R_{2r}(n,h)x^h q^n
	=
	\dfrac{1}{(xq;q^2)_{\infty}}
	\dfrac{(q^{2r};q^{2r})_{\infty}}
	{(q^2;q^2)_{\infty}}
	=
	\dfrac{(q^{2r};q^{2r})_{\infty}}
	{(xq;q^2)_{\infty}(q^2;q^2)_{\infty}},$$
	which completes the proof.
\end{proof}
\begin{remark}
	Putting $x=1$ in \eqref{RGF}, and comparing the coefficients of $q^n$, we obtain
	\begin{equation}\label{rmk1}
		\mathcal R_{2r}(n)
		=
		\sum_{h=0}^{\infty}\mathcal R_{2r}(n,h).
	\end{equation}
\end{remark}

\begin{definition}\label{def1}
Let $h\geq0$, $r\geq 2$, and $s\geq 1$ be integers. Let
$BELO_{r,s}(n,h)$  be the number of partitions of nonnegative integer $n$ in which every even part occurs at most $r-1$ times, the largest part is odd and occurs at least $s$ times, and the partition contains exactly $h$ odd parts, counted with multiplicity. 
\end{definition}
For example, $BELO_{3,2}(8,2)=1$ with the partition $3+3+2$ and $BELO_{2,1}(7,3)=3$ with the relevant partitions $5+1+1,\quad 3+3+1,\text{ and}\quad  3+2+1+1$.

\begin{theorem}
Let $r\geq2$, $s\geq1$, and $|xq|<1$ be integers. Then, we have
\begin{equation}\label{BGF}
	\sum_{n=0}^{\infty}\sum_{h=0}^{\infty}
	BELO_{r,s}(n,h)x^h q^n
	=
	\sum_{m=0}^{\infty}
	\dfrac{(q^{2r};q^{2r})_m(xq^{2m+1})^s}
	{(xq;q^2)_{m+1}(q^2;q^2)_m}.
\end{equation}   
\end{theorem}
\begin{proof}
Suppose that the largest part of a partition counted by
$BELO_{r,s}(n,h)$ is $2m+1$. Since the odd parts
$1,3,\ldots,2m-1$ can occur without restriction and the largest odd
part $2m+1$ occurs at least $s$ times, the contribution to $BELO_{r,s}(n,h)$ from the odd parts is given by 
\begin{equation}\label{BO1}
	\dfrac{(xq^{2m+1})^s}
	{(1-xq)(1-xq^3)\cdots(1-xq^{2m+1})}=
	\dfrac{(xq^{2m+1})^s}{(xq;q^2)_{m+1}},
\end{equation}
where the power of $x$ records the number of odd parts.

Since $2m+1$ is the largest part in the partition of $n$, the possible even parts are
$2,4,\ldots,2m$. Since each even part $2j$ can occur
$0,1,\ldots,r-1$ times, their contribution is given by
\begin{equation}\label{BE1}
	\prod_{j=1}^{m}
	\left(1+q^{2j}+q^{4j}+\cdots+q^{2(r-1)j}\right)
	=
	\prod_{j=1}^{m}\dfrac{1-q^{2rj}}{1-q^{2j}}
	=
	\dfrac{(q^{2r};q^{2r})_m}{(q^2;q^2)_m}.
\end{equation}
Multiplying \eqref{BO1} and \eqref{BE1}, and then summing over
$m\geq0$, we arrive at the desired result.
\end{proof}
\begin{remark}
Setting $x=1$ in \eqref{BGF} and \begin{equation}\label{rmrk1e1}
	\sum_{h=0}^{\infty}BELO_{r,s}(n,h)=BELO_{r,s}(n),
\end{equation} we obtain
\begin{equation}\label{gn1}
	\sum_{n=0}^{\infty}BELO_{r,s}(n)q^n=
	\sum_{m=0}^{\infty}
	\dfrac{(q^{2r};q^{2r})_m q^{s(2m+1)}}
	{(q;q^2)_{m+1}(q^2;q^2)_m}=
	\sum_{m=0}^{\infty}
	\dfrac{(q^{2r};q^{2r})_m q^{s(2m+1)}}
	{(q;q)_{2m+1}}.
\end{equation}
Note that,  $BELO_{r,s}(n)$ counts the number of partitions of the nonnegative integer $n$ in which every even part occurs at most $r-1$ times and the largest part is odd and occurs at least $s$ times.   
\end{remark}

\begin{definition}\label{def4}
	Let $h\geq0$, $r\geq2$, and $s\geq1$ be integers. For every integer
	$n\geq1$, let $BELO_{r,\leq s}(n,h)$ to be the number of
	partitions of $n$ in which every even part occurs at most $r-1$ times,
	the largest part is odd and occurs at most $s$ times, and the partition
	contains exactly $h$ odd parts, counted with multiplicity. For convenience, we assume  that 
	$BELO_{r,\leq s}(0,h)=0$.
\end{definition}
For example, $BELO_{3,\leq2}(8,2)=5$ with the relevant partitions
$7+1$, $5+3$, $5+2+1$, $3+3+2$ and $3+2+2+1$,  and
$BELO_{2,\leq1}(7,3)=2$ with the  partitions 
$5+1+1$ and $3+2+1+1$.
\begin{theorem}\label{th16}
	Let $r\geq 2$, $s\geq 1$, and $|xq|<1$ be integers, then 
	\begin{equation}\label{BL16}
		\sum_{n=0}^{\infty}\sum_{h=0}^{\infty}
		BELO_{r,\leq s}(n,h)x^{h}q^{n}
		=
		\sum_{m=0}^{\infty}
		\dfrac{
			(q^{2r};q^{2r})_{m}\,xq^{2m+1}
			\left(1-x^{s}q^{s(2m+1)}\right)
		}{
			(xq;q^{2})_{m+1}(q^{2};q^{2})_{m}
		}.
	\end{equation}
\end{theorem}
\begin{proof}
	The partitions counted by $BELO_{r,s+1}(n,h)$ are precisely those
	partitions counted by $BELO_{r,1}(n,h)$ in which the largest odd part
	occurs more than $s$ times. By Definitions \ref{def1} and \ref{def4}, for $n\geq0$, we obtain
\begin{equation}\label{c2}BELO_{r,\leq s}(n,h)=
	BELO_{r,1}(n,h)-BELO_{r,s+1}(n,h).\end{equation}
	Applying \eqref{BGF}  with $s$ replaced by  $1$
	and $s+1$ i the right hand side of \eqref{c2} , we  complete the proof.
\end{proof}
\begin{remark}
	Setting  $x=1$ and
	\begin{equation}\label{beloles}
		\sum_{h=0}^{\infty}BELO_{r,\leq s}(n,h)
		=
		BELO_{r,\leq s}(n)
	\end{equation}
	in \eqref{BL16}, we obtain
		$$\sum_{n=0}^{\infty}BELO_{r,\leq s}(n)q^{n}
		=
		\sum_{m=0}^{\infty}
		\dfrac{(q^{2r};q^{2r})_{m}q^{2m+1}
			\left(1-q^{s(2m+1)}\right)}
		{(q;q)_{2m+1}}.$$
	Note that,  $BELO_{r,\leq s}(n)$ counts the partitions of $n$ in which
	every even part occurs at most $r-1$ times, the largest part is odd,
	and the largest part occurs at most $s$ times.
\end{remark}

\begin{definition}
	Let $r\geq 2$ and $t\geq 1$ be integers. For integer $n\geq0$, let $BELO_{r,=t}(n)$ denote the number of partitions counted by $BELO_{r,t}(n)$ in which the largest odd part occurs exactly $t$ times. 
\end{definition}
For example, $BELO_{3,=2}(8)=2$ with relevant partitions 
$3+3+2$
and
$3+3+1+1,$
and $BELO_{2,=1}(7)=5$ with the relevant partitions given by 
$7, 
5+2, 
5+1+1, 
3+2+1+1,$ and 
$3+1+1+1+1.$
\begin{theorem}
	Let $r\geq2$ and $t\geq1$ be integers. Then,  we have 
	\begin{equation*}
		\sum_{n=0}^{\infty} BELO_{r, =t}(n)q^n=\sum_{m=0}^{\infty}
		\dfrac{(q^{2r};q^{2r})_m q^{t(2m+1)}}{(q;q)_{2m}}.
	\end{equation*}
\end{theorem}
\begin{proof}It is easy to see that
	\begin{equation}\label{G9a}
		BELO_{r,=t}(n)
		=
		BELO_{r,t}(n)
		-
		BELO_{r,t+1}(n).
	\end{equation}
Employing \eqref{gn1} with $s$ replaced by $t$ and $t+1$ in the right hand side of \eqref{G9a} and simplifying, we complete the proof.
\end{proof}

\section{Identities relating the  partition functions  and $2r$-regular partition function}
This section is devoted to proving identities relating our partition functions and the $2r$-regular partition function.

\begin{theorem} 
Let $h\geq0$ and $r\geq2$ be integers. Then for $n\geq1$, we have 
\begin{equation*}
\mathcal R_{2r}(n,h)
=
BELO_{r,1}(n,h)
+
\sum_{j=1}^{r-1}
BELO_{r,j}(n-j,h+j),
\end{equation*}
where $BELO_{r,s}(n, h)=0$ for any integer $n<0$.
\end{theorem}
\begin{proof}
From \eqref{RGF}, we note that
\begin{align}
\sum_{n=0}^{\infty}\sum_{h=0}^{\infty}
\mathcal R_{2r}(n,h)x^hq^n-1
&=\dfrac{(q^{2r};q^{2r})_\infty}
{(xq;q^2)_\infty(q^2;q^2)_\infty}-1
\notag\\
&=\sum_{m=0}^{\infty}
\left(
\dfrac{(q^{2r};q^{2r})_{m+1}}
{(xq;q^2)_{m+1}(q^2;q^2)_{m+1}}
-
\dfrac{(q^{2r};q^{2r})_m}
{(xq;q^2)_m(q^2;q^2)_m}
\right)
\notag\\
&=\sum_{m=0}^{\infty}
\dfrac{(q^{2r};q^{2r})_m}
{(xq;q^2)_m(q^2;q^2)_m}
\left(
\dfrac{1-q^{2r(m+1)}}
{(1-xq^{2m+1})(1-q^{2m+2})}-1
\right)
\notag\\
&=\sum_{m=0}^{\infty}
\dfrac{(q^{2r};q^{2r})_m}
{(xq;q^2)_m(q^2;q^2)_m}\times
\dfrac{xq^{2m+1}-xq^{4m+3}
+q^{2m+2}-q^{2r(m+1)}}
{(1-xq^{2m+1})(1-q^{2m+2})}
\notag\\
&=\sum_{m=0}^{\infty}
\dfrac{(q^{2r};q^{2r})_m}
{(xq;q^2)_m(q^2;q^2)_m}\times
\dfrac{(1-q^{2m+2})
\left(
xq^{2m+1}
+\displaystyle\sum_{j=1}^{r-1}q^{j(2m+2)}
\right)}
{(1-xq^{2m+1})(1-q^{2m+2})}
\notag\\
&=\sum_{m=0}^{\infty}
\dfrac{(q^{2r};q^{2r})_m}
  {(xq;q^2)_{m+1}(q^2;q^2)_m}
\left(
xq^{2m+1}
+\sum_{j=1}^{r-1}q^{j(2m+2)}
\right)
\notag\\
&=\sum_{m=0}^{\infty}
\dfrac{(q^{2r};q^{2r})_m(xq^{2m+1})}
{(xq;q^2)_{m+1}(q^2;q^2)_m}
+\sum_{j=1}^{r-1}q^jx^{-j}
\sum_{m=0}^{\infty}
\dfrac{(q^{2r};q^{2r})_m(xq^{2m+1})^j}
{(xq;q^2)_{m+1}(q^2;q^2)_m}
\notag\\
&\label{m1e1}=\sum_{n=0}^{\infty}\sum_{h=0}^{\infty}
BELO_{r,1}(n,h)x^hq^n
+\sum_{j=1}^{r-1}q^jx^{-j}
\sum_{n=0}^{\infty}\sum_{h=0}^{\infty}
BELO_{r,j}(n,h)x^hq^n.
\end{align}
Equating  the coefficients of $x^hq^n$ on both sides of \eqref{m1e1}, we complete the proof.
\end{proof}
\begin{corollary}\label{cor}
Let $r\geq2$ be any integer. Then, for any integer $n\geq1$, we have
\begin{equation}\label{BC3}
\mathcal R_{2r}(n)
=
BELO_{r,1}(n)
+
\sum_{j=1}^{r-1}
BELO_{r,j}(n-j).
\end{equation}
where $BELO_{r,s}(n)=0$ for any integer $n<0$.
\end{corollary}
\begin{proof}
Putting $x=1$ in \eqref{m1e1}  and then employing \eqref{rmk1} and \eqref{rmrk1e1}, we obtain 
\begin{align}
\dfrac{(q^{2r};q^{2r})_{\infty}}
{(q;q)_{\infty}}-1
&\label{BC2}=
\sum_{n=0}^{\infty}
BELO_{r,1}(n)q^n+
\sum_{j=1}^{r-1}q^j
\sum_{n=0}^{\infty}
BELO_{r,j}(n)q^n.
\end{align}
Equating the coefficients of $q^n$ on both sides of  \eqref{BC2}, we arroive at the desired result. 
\end{proof}

\begin{theorem}\label{t8}
Let $r\geq2$ and  $t\geq2$ be integers. Then for $n\geq0$,  we have
\begin{align}
\mathcal R_{2r}(n)
&=BELO_{r,1}(n)
+\sum_{j=1}^{t-1}BELO_{r,j}(n-j)+BELO_{r,t}(n+t)
-BELO_{r,t+1}(n+t)\notag\\
&\label{t12e1}\quad\hspace{2cm}-\sum_{j=0}^{t-1}BELO_{r,r+j}(n-r-j),
\end{align}
where $BELO_{r,s}(n)=0$ for any integer $n<0$.
\end{theorem}
\begin{proof}
From \eqref{gn1}, we obtain
\begin{align}
&\hspace{-2cm}q^{-t}\sum_{n=0}^{\infty}
\bigl(BELO_{r,t}(n)-BELO_{r,t+1}(n)\bigr)q^n
\notag\\
&\hspace{1cm}=
\sum_{m=0}^{\infty}
\dfrac{(q^{2r};q^{2r})_m q^{t(2m+1)-t}}
{(q;q)_{2m+1}}
-
\sum_{m=0}^{\infty}
\dfrac{(q^{2r};q^{2r})_m q^{(t+1)(2m+1)-t}}
{(q;q)_{2m+1}}
\notag\\
&\hspace{1cm}=
\sum_{m=0}^{\infty}
\dfrac{(q^{2r};q^{2r})_m q^{2mt}
\left(1-q^{2m+1}\right)}
{(q;q)_{2m}\left(1-q^{2m+1}\right)}
\notag\\
&\hspace{1cm}=
\sum_{m=0}^{\infty}
\dfrac{(q^{2r};q^{2r})_m q^{2mt}}
{(q;q)_{2m}}
\notag\\
&\hspace{1cm}=
1+
\sum_{m=0}^{\infty}
\dfrac{(q^{2r};q^{2r})_{m+1}
q^{2t(m+1)}}
{(q;q)_{2m+2}}
\notag\\
&\hspace{1cm}=
1+
\sum_{m=0}^{\infty}
\dfrac{(q^{2r};q^{2r})_m
q^{2t(m+1)}
\left(1-q^{2r(m+1)}\right)}
{(q;q)_{2m+1}\left(1-q^{2m+2}\right)}
\notag\\
&\hspace{1cm}=
1+
\sum_{m=0}^{\infty}
\dfrac{(q^{2r};q^{2r})_m q^{2t(m+1)}}
{(q;q)_{2m+1}}
\sum_{j=0}^{r-1}q^{2j(m+1)}
\notag\\
&\hspace{1cm}=
1+
\sum_{j=0}^{r-1}q^{t+j}
\sum_{m=0}^{\infty}
\dfrac{(q^{2r};q^{2r})_m
q^{(t+j)(2m+1)}}
{(q;q)_{2m+1}}
\notag\\
&\label{l7}\hspace{1cm}=
1+
\sum_{j=0}^{r-1}q^{t+j}
\sum_{n=0}^{\infty}BELO_{r,t+j}(n)q^n,
\end{align}
which implies
\begin{align}
&\hspace{.5cm}\sum_{n=0}^{\infty}BELO_{r,1}(n)q^{n}
+\sum_{j=1}^{t-1}q^{j}
 \sum_{n=0}^{\infty}BELO_{r,j}(n)q^{n}+
 q^{-t}\sum_{n=0}^{\infty}
 \left(
 BELO_{r,t}(n)-BELO_{r,t+1}(n)
 \right)q^{n}
\nonumber\\
&\hspace{.5cm}
-
\sum_{j=0}^{t-1}q^{r+j}
\sum_{n=0}^{\infty}BELO_{r,r+j}(n)q^{n}
=
1+\sum_{n=0}^{\infty}BELO_{r,1}(n)q^{n}
+\sum_{j=1}^{t-1}q^{j}
 \sum_{n=0}^{\infty}BELO_{r,j}(n)q^{n}
\nonumber\\
&\hspace{6cm}+
\sum_{j=t}^{t+r-1}q^{j}
\sum_{n=0}^{\infty}BELO_{r,j}(n)q^{n}
-
\sum_{j=r}^{r+t-1}q^{j}
\sum_{n=0}^{\infty}BELO_{r,j}(n)q^{n}
\nonumber\\
&\hspace{6cm}=
1+\sum_{n=0}^{\infty}BELO_{r,1}(n)q^{n}
+\sum_{j=1}^{r-1}q^{j}
 \sum_{n=0}^{\infty}BELO_{r,j}(n)q^{n}
\nonumber\\
&\label{l1}\hspace{6cm}=
\sum_{n=0}^{\infty}R_{2r}(n)q^{n},
\end{align}
where the last equality follows from \eqref{BC2}. Equating the coefficients of $q^n$ in \eqref{l1}, we complete the proof.
\end{proof}
\begin{theorem}
Let $r\geq2$ and $s\geq1$ be integers. Then for every integer $n\geq0$, we have 
\begin{align}
&\hspace{-1.5cm}BELO_{r,s+1}(n)
+BELO_{r,s+1}(n-1)
+BELO_{r,s}(n-2s-1)\notag\\
&\label{T10}=
BELO_{r,s+2}(n)
+BELO_{r,s}(n-1)
+BELO_{r,s+r}(n-2s-r-1),
\end{align}
where $BELO_{r,s}(n)=0$ for any integer $n<0$.
\end{theorem}
\begin{proof}
For $m\geq1$,  we note that
\begin{equation*}
    (1-q^{2m})(1-q^{2m+1})
\dfrac{(q^{2r};q^{2r})_m}{(q;q)_{2m+1}}
=
\dfrac{(q^{2r};q^{2r})_m}{(q;q)_{2m-1}}
=
(1-q^{2rm})
\dfrac{(q^{2r};q^{2r})_{m-1}}{(q;q)_{2m-1}}. 
\end{equation*}
Therefore, 
\begin{equation}
\sum_{m=1}^{\infty}
(1-q^{2m})(1-q^{2m+1})
\dfrac{(q^{2r};q^{2r})_m q^{2sm}}
{(q;q)_{2m+1}}=
\sum_{m=1}^{\infty}
(1-q^{2rm})
\dfrac{(q^{2r};q^{2r})_{m-1}q^{2sm}}
{(q;q)_{2m-1}}\label{l2}    
\end{equation}
Simplifying \eqref{l2}, we obtain
\begin{align}
&\sum_{m=0}^{\infty}
\dfrac{(q^{2r};q^{2r})_m q^{2sm}}
{(q;q)_{2m+1}}
-(1+q)\sum_{m=0}^{\infty}
\dfrac{(q^{2r};q^{2r})_m q^{2(s+1)m}}
{(q;q)_{2m+1}}+q\sum_{m=0}^{\infty}
\dfrac{(q^{2r};q^{2r})_m q^{2(s+2)m}}
{(q;q)_{2m+1}}\notag\\
&\label{t10e1}=
q^{2s}\sum_{m=0}^{\infty}
\dfrac{(q^{2r};q^{2r})_m q^{2sm}}
{(q;q)_{2m+1}}
-q^{2s+2r}\sum_{m=0}^{\infty}
\dfrac{(q^{2r};q^{2r})_m q^{2(s+r)m}}
{(q;q)_{2m+1}}.
\end{align}
Employing \eqref{gn1} in \eqref{t10e1} and simplifying, we obtain
\begin{align}
	&\hspace{-1cm}(1+q)\sum_{n=0}^{\infty}BELO_{r,s+1}(n)q^{n}
	+q^{2s+1}\sum_{n=0}^{\infty}BELO_{r,s}(n)q^{n} \notag\\
	&\label{T10E24}
	=\sum_{n=0}^{\infty}BELO_{r,s+2}(n)q^{n}
	+q\sum_{n=0}^{\infty}BELO_{r,s}(n)q^{n}
	+q^{2s+r+1}\sum_{n=0}^{\infty}
	BELO_{r,s+r}(n)q^{n}.
\end{align}
\eqref{T10E24} implies
\begin{align}
 &\hspace{-.5cm}\sum_{n=0}^{\infty}
 \left(
 BELO_{r,s+1}(n)
 +BELO_{r,s+1}(n-1)
 +BELO_{r,s}(n-2s-1)
 \right)q^{n} \notag\\
 &\label{T10E25}
 =\sum_{n=0}^{\infty}
 \left(
 BELO_{r,s+2}(n)
 +BELO_{r,s}(n-1)
 +BELO_{r,s+r}(n-2s-r-1)
 \right)q^{n}.
\end{align}
Equating  the coefficients of $q^{n}$ on both sides of \eqref{T10E25}, we arrive at the desired result.
\end{proof}
\begin{theorem}
Let $h\geq0$ and $r\geq2$ be integers. Then for any integer $n\geq1$,  we have
\begin{align}
\mathcal R_{2r}(n,h)
={}&BELO_{r,\leq h}(n,h)+\sum_{j=1}^{r-1}
\left(
BELO_{r,\leq h+j}(n-j,h+j)
-
BELO_{r,\leq j-1}(n-j,h+j)
\right),\notag
\end{align}
where $BELO_{r,\leq0}(n,h)=0$, and $BELO_{r,\leq s}(n,h)=0$, for any integer $n<0$.
\end{theorem}
\begin{proof}
For $t\geq1$, subtracting the case $s=t-1$ from the case $s=t$
of \eqref{BL16} and noting $BELO_{r,\leq0}(n,h)=0$, we obtain
\begin{align}
\sum_{n=0}^{\infty}\sum_{h=0}^{\infty}
\left(
BELO_{r,\leq t}(n,h)
-
BELO_{r,\leq t-1}(n,h)
\right)x^{h}q^{n}
&=
\sum_{m=0}^{\infty}
\dfrac{(q^{2r};q^{2r})_{m}
(xq^{2m+1})^{t}
(1-xq^{2m+1})}
{(xq;q^{2})_{m+1}(q^{2};q^{2})_{m}}
\notag\\
&\label{BL6}=
\sum_{m=0}^{\infty}
\dfrac{(q^{2r};q^{2r})_{m}
(xq^{2m+1})^{t}}
{(xq;q^{2})_{m}(q^{2};q^{2})_{m}}.
\end{align}  
Summing \eqref{BL6} over $t\geq k\ge1$ and
using \eqref{BGF}, we obtain
\begin{align}
\sum_{t=k}^{\infty}
\sum_{n=0}^{\infty}\sum_{h=0}^{\infty}
\left(
BELO_{r,\leq t}(n,h)
-
BELO_{r,\leq t-1}(n,h)
\right)x^{h}q^{n}
&=
\sum_{m=0}^{\infty}
\dfrac{(q^{2r};q^{2r})_{m}}
{(xq;q^{2})_{m}(q^{2};q^{2})_{m}}
\sum_{t=k}^{\infty}(xq^{2m+1})^{t}
\notag\\
&=
\sum_{m=0}^{\infty}
\dfrac{(q^{2r};q^{2r})_{m}
(xq^{2m+1})^{k}}
{(1-xq^{2m+1})(xq;q^{2})_{m}(q^{2};q^{2})_{m}}
\notag\\
&=
\sum_{m=0}^{\infty}
\dfrac{(q^{2r};q^{2r})_{m}
(xq^{2m+1})^{k}}
{(xq;q^{2})_{m+1}(q^{2};q^{2})_{m}}
\notag\\
&\label{BL7}=
\sum_{n=0}^{\infty}\sum_{h=0}^{\infty}
BELO_{r,k}(n,h)x^{h}q^{n}.
\end{align}
Substituting \eqref{BL7} with cases $k=1$ and $k=j$ in \eqref{m1e1}, we obtain
\begin{align}
\sum_{n=0}^{\infty}\sum_{h=0}^{\infty}
\mathcal R_{2r}(n,h)x^{h}q^{n}-1
&=
\sum_{t=1}^{\infty}
\sum_{n=0}^{\infty}\sum_{h=0}^{\infty}
\left(
BELO_{r,\leq t}(n,h)
-
BELO_{r,\leq t-1}(n,h)
\right)x^{h}q^{n}
\notag\\
&\quad+
\sum_{j=1}^{r-1}q^{j}x^{-j}
\sum_{t=j}^{\infty}
\sum_{n=0}^{\infty}\sum_{h=0}^{\infty}
\left(
BELO_{r,\leq t}(n,h)
-
BELO_{r,\leq t-1}(n,h)
\right)x^{h}q^{n}
\notag\\
&=
\sum_{n=0}^{\infty}\sum_{h=0}^{\infty}
\Bigg(
\sum_{t=1}^{h}
\left(
BELO_{r,\leq t}(n,h)
-
BELO_{r,\leq t-1}(n,h)
\right)
\notag\\
&\label{t18e1}\quad+
\sum_{j=1}^{r-1}\sum_{t=j}^{h+j}
\left(
BELO_{r,\leq t}(n-j,h+j)
-
BELO_{r,\leq t-1}(n-j,h+j)
\right)
\Bigg)x^{h}q^{n}.
\end{align}
Comparing the coefficients of $x^{h}q^{n}$  on both sides of \eqref{t18e1}, we obtain
\begin{align}
\mathcal R_{2r}(n,h)
&=
\sum_{t=1}^{h}
\left(
BELO_{r,\leq t}(n,h)
-
BELO_{r,\leq t-1}(n,h)
\right)
\notag\\
&\quad+
\sum_{j=1}^{r-1}\sum_{t=j}^{h+j}
\left(
BELO_{r,\leq t}(n-j,h+j)
-
BELO_{r,\leq t-1}(n-j,h+j)
\right)
\notag\\
&=
BELO_{r,\leq h}(n,h)
+
\sum_{j=1}^{r-1}
\left(
BELO_{r,\leq h+j}(n-j,h+j)
-
BELO_{r,\leq j-1}(n-j,h+j)
\right).
\notag
\end{align}
Hence, the proof is complete.
\end{proof}
\begin{corollary}
Let $r\geq2$ be integer. Then for any integer $n\geq1$, we have
$$\mathcal R_{2r}(n)
=BELO_{r,\leq n}(n)+\sum_{j=1}^{\min\{r-1,n-1\}}
\left(
BELO_{r,\leq n-j}(n-j)
-
BELO_{r,\leq j-1}(n-j)
\right).$$
where $BELO_{r,\leq0}(n)=0$, and $BELO_{r,\leq s}(n)=0$ for any  integer $n<0$.
\end{corollary}
\begin{proof}
Setting $x=1$ in \eqref{t18e1} and then employing \eqref{beloles}, we obtain
\begin{align}
\sum_{n=0}^{\infty}R_{2r}(n)q^{n}-1
={}&
\sum_{t=1}^{\infty}\sum_{n=0}^{\infty}
\left(
BELO_{r,\leq t}(n)
-
BELO_{r,\leq t-1}(n)
\right)q^{n}
\nonumber\\
&\label{BL9}+
\sum_{j=1}^{r-1}q^{j}
\sum_{t=j}^{\infty}\sum_{n=0}^{\infty}
\left(
BELO_{r,\leq t}(n)
-
BELO_{r,\leq t-1}(n)
\right)q^{n}.
\end{align}
Equating the coefficients of $q^{n}$ on both sides of  \eqref{BL9}, we obtain
\begin{align}
R_{2r}(n)
=&
\sum_{t=1}^{\infty}
\left(
BELO_{r,\leq t}(n)
-
BELO_{r,\leq t-1}(n)
\right)
\nonumber\\
&\label{BL10}+
\sum_{j=1}^{\min\{r-1,n-1\}}
\sum_{t=j}^{\infty}
\left(
BELO_{r,\leq t}(n-j)
-
BELO_{r,\leq t-1}(n-j)
\right).
\end{align}
Since $
BELO_{r,\leq t}(n)-BELO_{r,\leq t-1}(n)
$ counts the number of partitions of any positive integer $n$ in which every even part appears at most $r-1$ times and  the largest odd part occurs exactly $t$ times,  the  multiplicity of the   largest odd part in such partitions cannot exceed $n$.  Therefore, 
\begin{align}
\sum_{t=1}^{\infty}
\left(
BELO_{r,\leq t}(n)
-
BELO_{r,\leq t-1}(n)
\right)
&\label{c19e1}=BELO_{r,\leq n}(n)\\\intertext{and}
\sum_{t=j}^{\infty}
\left(
BELO_{r,\leq t}(n-j)
-
BELO_{r,\leq t-1}(n-j)
\right)
&\label{c19e2}=
BELO_{r,\leq n-j}(n-j)
-
BELO_{r,\leq j-1}(n-j).
\end{align}
Employing \eqref{c19e1} and \eqref{c19e2} in  \eqref{BL10}, we complete the proof.
\end{proof}

\begin{theorem}
Let $r\geq2$ and $s\geq1$ be integers. Then for any integer $n\geq0$, we have 
\begin{align}
&\hspace{-.5cm}BELO_{r,\leq s}(n)
+BELO_{r,\leq s}(n-1)
+BELO_{r,\leq s-1}(n-2s-1)\notag\\
&\hspace{1cm}\quad=
BELO_{r,\leq s+1}(n)
+BELO_{r,\leq s-1}(n-1)
+BELO_{r,\leq s+r-1}(n-2s-r-1)\notag\\
&
\qquad\hspace{1cm}+BELO_{r,\leq n-2s-1}(n-2s-1)
-BELO_{r,\leq n-2s-r-1}(n-2s-r-1),\notag
\end{align}
where $BELO_{r,\leq0}(n)=0$ and $BELO_{r,\leq s}(n)=0$ for any integer $n<0$.
\end{theorem}
\begin{proof}
Note that,
\begin{align}
BELO_{r,k}(n)
&=
\sum_{h=0}^{\infty}BELO_{r,k}(n,h)
\notag\\
&=
\sum_{h=0}^{\infty}
\left(
BELO_{r,\leq n}(n,h)
-
BELO_{r,\leq k-1}(n,h)
\right)
\notag\\
&\label{C19B}=
BELO_{r,\leq n}(n)-BELO_{r,\leq k-1}(n).
\end{align}
Using  \eqref{C19B}  in
\eqref{T10} with the corresponding values of $n$ and $k$, we obtain
\begin{align}
&  \left(
BELO_{r,\leq n}(n)-BELO_{r,\leq s}(n)
\right)
+
\left(
BELO_{r,\leq n-1}(n-1)-BELO_{r,\leq s}(n-1)
\right)
\notag\\
&\hspace{1.7cm}\qquad+
\left(
BELO_{r,\leq n-2s-1}(n-2s-1)
-
BELO_{r,\leq s-1}(n-2s-1)
\right)
\notag\\
&\hspace{1.4cm}\quad=
\left(
BELO_{r,\leq n}(n)-BELO_{r,\leq s+1}(n)
\right)
+
\left(
BELO_{r,\leq n-1}(n-1)
-
BELO_{r,\leq s-1}(n-1)
\right)
\notag\\
&\label{l4}\qquad \hspace{1.7cm}+
\left(
BELO_{r,\leq n-2s-r-1}(n-2s-r-1)
-
BELO_{r,\leq s+r-1}(n-2s-r-1)
\right).
\end{align}
Simplifying \eqref{l4}, we  arrive at the desired result.
\end{proof}
\begin{theorem}
Let $r\geq2$ and $t\geq2$ be integers. Then for any integer $n\geq0$, we have
\begin{align}
\mathcal R_{2r}(n)
=BELO_{r,\leq n}(n)
&+\sum_{j=1}^{t-1}
\left(
BELO_{r,\leq n-j}(n-j)
-
BELO_{r,\leq j-1}(n-j)
\right)\notag\\
&+BELO_{r,\leq t}(n+t)
-BELO_{r,\leq t-1}(n+t)\notag\\
&-\sum_{j=0}^{t-1}
\left(
BELO_{r,\leq n-r-j}(n-r-j)
-
BELO_{r,\leq r+j-1}(n-r-j)
\right),\notag
\end{align}
where $BELO_{r,\leq0}(n)=0$, and $BELO_{r,\leq s}(n)=0$ for any  integer $n<0$.
\end{theorem}
\begin{proof}
Employing \eqref{C19B} in \eqref{t12e1}, we easiy arrive at the desired result.
\end{proof}

\begin{theorem}
	Let $r\geq2$ and $n\geq1$ be integers. Then
		$$\mathcal{R}_{2r}(n)
		=
		\sum_{t=1}^{\infty}
		\sum_{j=0}^{\min\{t,r-1\}}
		BELO_{r,=t}(n-j),$$
	where $BELO_{r,=t}(n)=0$ for any  integer $n<0$.
\end{theorem}
\begin{proof}
	From \eqref{BC3}, we note that
	\begin{align}
		\mathcal{R}_{2r}(n)
		&=
		BELO_{r,1}(n)
		+
		\sum_{j=1}^{r-1}
		BELO_{r,j}(n-j)
		\notag\\
		&=
		\sum_{t=1}^{\infty}
		BELO_{r,=t}(n)
		+
		\sum_{j=1}^{r-1}
		\sum_{t=j}^{\infty}
		BELO_{r,=t}(n-j)
		\notag\\
		&=
		\sum_{t=1}^{\infty}
		BELO_{r,=t}(n)
		+
		\sum_{t=1}^{\infty}
		\sum_{j=1}^{\min\{t,r-1\}}
		BELO_{r,=t}(n-j)
		\notag	\end{align}
		\begin{align}
		&~=
		\sum_{t=1}^{\infty}
		\left(
		BELO_{r,=t}(n)
		+
		\sum_{j=1}^{\min\{t,r-1\}}
		BELO_{r,=t}(n-j)
		\right)
		\notag\\
		&~=
		\sum_{t=1}^{\infty}
		\sum_{j=0}^{\min\{t,r-1\}}
		BELO_{r,=t}(n-j).
		\notag
	\end{align}
Hence, the proof is complete.
\end{proof}

\begin{theorem}
Let $r\geq2$ and $t\geq1$ be integers. Then for any integer
$n\geq0$, we have
\begin{align}
&BELO_{r,=t+1}(n-1)+BELO_{r,=t+1}(n-2)
+BELO_{r,=t}(n-2t-3)
\notag\\
&\qquad=
BELO_{r,=t+2}(n)+BELO_{r,=t}(n-3)
+BELO_{r,=t+r}(n-2t-r-3),\notag
\end{align}
where $BELO_{r,=t}(k)=0$, if $k$ is a negative integer.
\end{theorem}
\begin{proof}
Setting $s=t$ in \eqref{T10E24}, we obtain
\begin{align}
&(1+q)\sum_{n=0}^{\infty}BELO_{r,t+1}(n)q^n
+q^{2t+1}\sum_{n=0}^{\infty}BELO_{r,t}(n)q^n\notag\\
&\label{t13e1}=
\sum_{n=0}^{\infty}BELO_{r,t+2}(n)q^n
+q\sum_{n=0}^{\infty}BELO_{r,t}(n)q^n
+q^{2t+r+1}
\sum_{n=0}^{\infty}BELO_{r,t+r}(n)q^n.
\intertext{Similarly, taking $s=t+1$ in \eqref{T10E24}, we obtain}
&(1+q)\sum_{n=0}^{\infty}BELO_{r,t+2}(n)q^n
+q^{2t+3}\sum_{n=0}^{\infty}BELO_{r,t+1}(n)q^n
\notag\\
&\label{t13e2}=
\sum_{n=0}^{\infty}BELO_{r,t+3}(n)q^n
+q\sum_{n=0}^{\infty}BELO_{r,t+1}(n)q^n
+q^{2t+r+3}
\sum_{n=0}^{\infty}BELO_{r,t+r+1}(n)q^n.
\intertext{Multiplying  \eqref{t13e1} by $q^2$, we obtain }
&(q^2+q^3)
\sum_{n=0}^{\infty}BELO_{r,t+1}(n)q^n
+q^{2t+3}
\sum_{n=0}^{\infty}BELO_{r,t}(n)q^n
\notag\\
&\label{t13e3}=
q^2\sum_{n=0}^{\infty}BELO_{r,t+2}(n)q^n
+q^3\sum_{n=0}^{\infty}BELO_{r,t}(n)q^n
+q^{2t+r+3}
\sum_{n=0}^{\infty}BELO_{r,t+r}(n)q^n.
\intertext{Subtracting \eqref{t13e2} from \eqref{t13e3}, we obtain}
&(q^2+q^3)
\sum_{n=0}^{\infty}BELO_{r,t+1}(n)q^n
+q^{2t+3}
\sum_{n=0}^{\infty}BELO_{r,t}(n)q^n
-(1+q)
\sum_{n=0}^{\infty}BELO_{r,t+2}(n)q^n\notag\\
&\qquad
-q^{2t+3}
\sum_{n=0}^{\infty}BELO_{r,t+1}(n)q^n
\notag\\
&=
q^2\sum_{n=0}^{\infty}BELO_{r,t+2}(n)q^n
+q^3\sum_{n=0}^{\infty}BELO_{r,t}(n)q^n
\notag\\
&\quad
+q^{2t+r+3}
\sum_{n=0}^{\infty}BELO_{r,t+r}(n)q^n
-\sum_{n=0}^{\infty}BELO_{r,t+3}(n)q^n
\notag\\
&\label{l5}\quad
-q\sum_{n=0}^{\infty}BELO_{r,t+1}(n)q^n
-q^{2t+r+3}
\sum_{n=0}^{\infty}BELO_{r,t+r+1}(n)q^n,
\\
\intertext{Rearranging the terms in \eqref{l5} and using \eqref{G9a}, we obtain}
&(q+q^2)
\sum_{n=0}^{\infty}BELO_{r,=t+1}(n)q^n
+q^{2t+3}
\sum_{n=0}^{\infty}BELO_{r,=t}(n)q^n
\notag\\
&=
\sum_{n=0}^{\infty}BELO_{r,=t+2}(n)q^n
+q^3\sum_{n=0}^{\infty}BELO_{r,=t}(n)q^n
+q^{2t+r+3}
\sum_{n=0}^{\infty}BELO_{r,=t+r}(n)q^n
\notag\\
&\label{l6}=
\sum_{n=0}^{\infty}
\Bigl(
BELO_{r,=t+2}(n)
+BELO_{r,=t}(n-3)
+BELO_{r,=t+r}(n-2t-r-3)
\Bigr)q^n,
\end{align}
Equating the coefficients of $q^n$ on both sides \eqref{l6}, we arrive at the desired result.
\end{proof}

\begin{theorem}\label{tm}
Let $r\geq2$ and $s,t\geq1$ be integers. Then for  integer $n\geq0$, we have 
\begin{align}
BELO_{r,=t}(n)
={}&
\delta_{n,t}
+\sum_{j=0}^{r-1}
\Bigl(
BELO_{r,s}(n-2t-j)
-BELO_{r,\leq t+j-1}(n-2t-j)\notag\\
&\hspace{38mm}
+BELO_{r,\leq s-1}(n-2t-j)
\Bigr),\notag
\end{align}
   where $\delta_{n,t}$ is defined by
$$\delta_{n,t}
=
\begin{cases}
1, & \text{if } n=t,\\
0, & \text{if } n\neq t.
\end{cases}$$
Furthermore,  $BELO_{r,\leq0}(n)=0$, and $BELO_{r,\leq s}(n)=0$ for any integer $n<0$.
\end{theorem}
\begin{proof}
Employing \eqref{G9a} in \eqref{l7}, we obtain

\begin{equation}\label{BT3}
q^{-t}\sum_{n=0}^{\infty}BELO_{r,=t}(n)q^{n}
=
1+\sum_{j=0}^{r-1}q^{t+j}
\sum_{n=0}^{\infty}BELO_{r,t+j}(n)q^{n}.
\end{equation}
Substracting the case  $k=t+j$ from the case $k=s$ of  \eqref{C19B}
and simplifying, we obtain
\begin{align}
BELO_{r,t+j}(n)
={}&\label{BT7}
BELO_{r,s}(n)
-
BELO_{r,\leq t+j-1}(n)
+
BELO_{r,\leq s-1}(n).
\end{align}
Employing  \eqref{BT7} in \eqref{BT3} and then multiplying  by $q^{t}$, we obtain
\begin{align}
\sum_{n=0}^{\infty}BELO_{r,=t}(n)q^{n}
={}&q^{t}
+\sum_{j=0}^{r-1}q^{2t+j}
\sum_{n=0}^{\infty}
\Big(
BELO_{r,s}(n)
-BELO_{r,\leq t+j-1}(n)
\notag\\
&\label{BT2}\hspace{39mm}
+BELO_{r,\leq s-1}(n)
\Big)q^{n}.
\end{align}
Equating the coefficients of $q^n$ on both sides of \eqref{BT2}, we arrive at the desired result.
\end{proof}
\begin{remark}\label{rem}
Noting that,  $\operatorname{DE1}(n)=BELO_{2,1}(n)$, 
$\operatorname{DE2}(n)=BELO_{2,2}(n)$, and $\operatorname{DE3}(n)=BELO_{2,=1}(n)$,  Corollary 1, 2, 3 and 4  due to Andrews and El Bachraoui \cite{DE} stated in the introduction section follow from Corollary \ref{cor} and Theorem \ref{tm} as particular cases. Precisely,  
Corollary 1, 2  and 3  of \cite{DE} follow  from Corolarry \ref{cor} and Corollary 4 of \cite{DE}  follows from Theorem \ref{tm} proved in Section 3 of this paper. \end{remark}

\section*{\bf Declarations}
\noindent\textbf{Funding}: This research did not receive funding.\\
\noindent{\bf Author Contributions.} Both authors contributed equally to this work.

\noindent{\bf Conflict of Interest.} The authors declare that there is no conflict of interest regarding the publication of this paper.
 
\noindent{\bf Human and Animal Rights.} The authors declare that there is no research involving human participants or animals in the context of this paper.	
 
\noindent{\bf Data Availability Statement.} Data sharing is not applicable to this paper as no datasets were generated or analyzed during the current study.

\end{document}